\documentclass{elsart}

\usepackage{amsmath,amssymb}
\begin{document}

\begin{frontmatter}
\hsize=6.5in


\title{A generalization of the Widder potential transform and applications}

\author[label1]{Ne\c se Dernek},
\author[label2]{Veli Kurt},
\author[label2]{Y\i lmaz \c{S}im\c{s}ek},
\author[label3]{Osman Y\"urekli \corauthref{cor1}}

\address[label1]{Department of Mathematics, University of Marmara, TR-34722, Kadik\"oy, Istanbul, Turkey}
\address[label2]{Department of Mathematics, Akdeniz University, Antalya, Turkey}
\address[label3]{Department of Mathematics, Ithaca College, Ithaca, NY 14850, USA}

\corauth[cor1]{Corresponding author}

\begin{abstract}

In the present paper the authors consider the $\mathcal{P}_{\nu,2}$-transform as a generalization of the Widder potential transform and the Glasser transform. The $\mathcal{P}_{\nu,2}$-transform is obtained as an iteration of the the $\mathcal{L}_{2}$-transform with itself. Many identities involving these transforms are given. By making use of these identities, a number of new Parseval-Goldstein type identities are obtained for these and many other well-known integral transforms. The identities proven in this paper are shown to give rise to useful corollaries for evaluating infinite integrals of special functions. Some examples are also given as illustration of the results presented here.

\end{abstract}

\begin{keyword}
Laplace transforms\sep $\mathcal{L}_{2}$-transforms\sep Widder potential transforms \sep Glasser transforms \sep $\mathcal{P}_{\nu,2}$-transforms \sep Hankel transforms \sep $\mathcal{K}_{\nu}$-transforms \sep Parseval-Goldstein type theorems.
\par\bigskip
{\bf 2000 Mathematics Subject Classification.} Primary 44A10\sep 44A15; Secondary 33C10\sep 44A35.
\end{keyword}

\end{frontmatter}

\section{Introduction}
\numberwithin{equation}{section}
Over a  decade ago, Sadek and Y\"urekli \cite{YS} presented a systematic account
of so-called the $\mathcal{L}_{2}$-transform:
\begin{equation}
\label{l2}
\mathcal{L}_{2}\big\{f(x);y\big\}
=\int_0^\infty x\,\exp\big(-x^{2}\,y^{2}\big)\,f(x)\,dx
\end{equation}
The $\mathcal{L}_{2}$-transform is related to the classical Laplace transform 
\begin{equation}
\label{lap}
\mathcal{L} \big\{f(x);y \big\}
=
\int_0^\infty  \exp(-x \,y )\,f(x)\,dx
\end{equation}
by means of the following relationships:
\begin{align}
\label{l2:l}
\mathcal{L}_{2}\left\{f(x);y\right\}
&=
\frac{1}{2}\,\mathcal{L}\Big\{f\big(\sqrt{x}\big);y^{2}\Big\},
\\
\label{ll2}
\mathcal{L}\left\{f(x);y\right\}
&=
2\,\mathcal{L}_{2}\Big\{f\big(x^{2}\big);\sqrt{y}\Big\}.
\end{align}
Subsequently, various Parseval-Goldstein type identities were given in (for example)  \cite{BDY}, \cite{DSY}, \cite{GL}, \cite{Y99a}, and \cite{Y99b} for the $\mathcal{L}_2$-transform. New solutions techniques were obtained for the Bessel differential equation in \cite{WY02} and the Hermite differential equation in \cite{WY03} using this integral transform. There are numerous analogous results in the literature on various integral transforms (see, for instance \cite{SY95}, \cite{Y89}, \cite{YG}, and \cite{YS98}).

Over four decades ago, Widder \cite{WD} presented a systematic account of the so-called Widder potential transform:
\begin{equation}
\label{w}
\mathcal{P}\big\{f(x);y\big\}
=
\int_{0}^{\infty}
\frac{x\,f(x)}{x^{2}+y^{2}}\,dx,
\end{equation}
which, by an exponential change of variables, becomes a convolution transform with kernel belonging to a general class investigated by Hirschman and Widder \cite{HW}.

Over three decades ago, Glasser \cite{GL} considered so-called the Glasser transform
\begin{equation}
\label{gl}
\mathcal{G}\big\{f(x);y\big\}
=
\int_{0}^{\infty}
\frac{f(x)}{\sqrt{x^{2}+y^{2}}}\,dx.
\end{equation}
Glasser gave the following Parseval-Goldstein type theorem (cf. \cite[p. 171, Eq. (4)]{GL})
\begin{equation}
\label{pggl}
\int_{0}^{\infty}f(x)\,\mathcal{G}\big\{g(y);x\big\}\,dx
=
\int_{0}^{\infty}g(x)\,\mathcal{G}\big\{f(y);x\big\}\,dx,
\end{equation}
and evaluated a number of infinite integrals involving Bessel functions. Additional results about the Glasser transform
can be found in Srivastava and Y\"urekli \cite{SY95} and  Kahramaner {\it et al.} \cite{KSY}.

In this article, we introduce a new generalization of the Widder potential transform and the Glasser transform. We establish potentially useful identities for so-called the $\mathcal{P}_{\nu,2}$-transform and several other known integral transforms. First of all, the $\mathcal{P}_{\nu,2}$-transform is defined by
\begin{equation}
\label{pn2}
\mathcal{P}_{\nu,2}\big\{f(x);y\big\}
=
\int_{0}^{\infty}
\frac{x\,f(x)}{\big(x^{2}+y^{2}\big)^{\nu}}\,dx.
\end{equation}
If we put $\nu=1$ in the definition \eqref{pn2} above, we obtain the Widder potential transform \eqref{w}:
\begin{equation}
\label{pn2:w}
\mathcal{P}_{1,2}\big\{f(x);y\big\}
=
\mathcal{P}\big\{f(x);y\big\}
=
\int_{0}^{\infty}
\frac{x\,f(x)}{x^{2}+y^{2}}\,dx.
\end{equation}
If we put $\nu=1/2$ in \eqref{pn2}, we obtain
the Glasser transform \eqref{gl}:
\begin{equation}  \label{pn2:g}
\mathcal{P}_{\frac{1}{2},2}\bigg\{\frac{f(x)}{x};y\bigg\} = \mathcal{G}\big\{f(x);y\big\} =
\int_{0}^{\infty} \frac{f(x)}{\sqrt{x^{2}+y^{2}}}\,dx.
\end{equation}
  The Hankel transform is defined by
\begin{equation}
\label{hankel}
\mathcal{H}_{\nu}\big\{f(x);y\big\}
=
\int_{0}^{\infty}
\sqrt{xy}\,J_{\nu}(x\,y)\,f(x)\,dx
\end{equation}
where $J_{\nu}(x)$ is the Bessel function of the first kind of order $\nu$, and the $\mathcal{K}$-transform is defined by
\begin{equation}
\label{kn}
\mathcal{K}_{\nu}\big\{f(x);y\big\}
=
\int_{0}^{\infty}
\sqrt{xy}\,K_{\nu}(x\,y)\,f(x)\,dx
\end{equation}
where $K_{\nu}$ is the Bessel function of the second kind of order $\nu$.

In Section 2 of this paper, we show that an iteration of $\mathcal{L}_{2}$-transform \eqref{l2} with itself is $\mathcal{P}_{\nu,2}$-transform \eqref{pn2}. Using this iteration identity, a number of new Parseval-Goldstein type identities are then obtained for these and many other well-known integral transforms. Our main theorem is shown to yield new identities for the integral transforms introduced above. As applications of the identities and the Theorem, some illustrative examples are also given.


\section{The Main Theorem}
\numberwithin{equation}{section}
In the following lemma, we give an iteration identity involving the $\mathcal{L}_{2}$-transform \eqref{l2} and the $\mathcal{P}_{\nu,2}$-transform \eqref{pn2}.
\begin{lem}\label{l1} 
The identity 
\begin{equation}
\label{l1:1}
\mathcal{L}_{2}\Big\{u^{2\nu-2}\,\mathcal{L}_{2}\Big\{g(x);u\Big\};y\Big\}
=
\frac{\Gamma(\nu)}{2}
\,\mathcal{P}_{\nu,2}\big\{g(x);y\big\},
\end{equation}
holds true, provided that $\Re(\nu)>0$ and the integrals involved converge absolutely.
\end{lem}
\begin{pf} 
Using the definition \eqref{l2} of the $\mathcal{L}_{2}$-transform,
we have
\begin{align}
\notag
\mathcal{L}_{2}\Big\{u^{2\nu-2}\,&\mathcal{L}_{2}\big\{g(x);u\big\};y\Big\}
\\&=
\label{l1:p1} 
\int_{0}^{\infty} 
u^{2\nu-1}\,\exp\big(-y^{2}\,u^{2}\big)\,
\bigg[
\int_{0}^{\infty} 
x\,\exp\big(-x^{2}\,u^{2}\big)\,g(x)\,dx
\bigg]\,du.
\end{align}
Changing the order of integration, which is permissible by absolute convergence of the integrals involved, and then using the definition \eqref{l2} of the $\mathcal{L}_{2}$-transform once more, we find from \eqref{l1:p1} that
\begin{align}
\notag
\mathcal{L}_{2}\Big\{u^{2\nu-2}
\,\mathcal{L}_{2}&\big\{g(x);u\big\};y\Big\}
\\
\notag
&=
\int_{0}^{\infty} 
x\,g(x)\,
\bigg[
\int_{0}^{\infty} 
u^{2\nu-1}\,\exp\Big[\big(-y^{2}+x^{2}\big)u^{2}\Big]\,du\bigg]\,dx
\\
\label{l1:p2}
&=
\int_{0}^{\infty} 
x\,g(x)\,
\mathcal{L}_{2}\Big\{u^{2\nu-2};\big(x^{2}+y^{2}\big)^{1/2}\Big\}
\,dx
\end{align}
Furthermore, we have
\begin{equation}
\label{l1:p3}
\mathcal{L}_{2}\Big\{u^{2\nu-2};\big(x^{2}+y^{2}\big)^{1/2}\Big\}=\frac{1}{2}\,\frac{\Gamma(\nu)}{\big(x^{2}+y^{2}\big)^{\nu}},
\end{equation}
where $\Re(\nu)>0$. Now the assertion \eqref{l1:1} follows from \eqref{l1:p2}, \eqref{l1:p3}, and the definition \eqref{pn2} of the $\mathcal{P}_{\nu,2}$-transform.
\qed
\end{pf}

Setting $\nu=1$ in \eqref{l1} of our Lemma \ref{l1} and using the relation \eqref{pn2:w}, we obtain the known identity (cf. \cite[p. 518, Eq. (2.1)]{YS}) contained in 
\begin{cor}\label{l2l2w}
We have\begin{equation}
\label{l2l2w:1}
\mathcal{L}_{2}\Big\{\mathcal{L}_{2}\big\{f(x);u\big\};y\Big\}
=
\mathcal{P}\big\{f(x);y\big\}.
\end{equation}
\end{cor}

In the following corollary we evaluate the $\mathcal{P}_{\nu,2}$-transform of the Bessel function of the first kind as an illustration of our Lemma \ref{l1}:
\begin{cor}\label{besselkMuNu}
We have
\begin{equation}
\label{besselkMuNuone}
\mathcal{P}_{\nu,2}\big\{x^{\mu}\,J_{\mu}(z\,x);y\big\}
=
\frac{1}{\Gamma(\nu)}\,\Big(\frac{z}{2}\Big)^{\nu-1}\,y^{\mu-\nu+1}\,K_{\nu-\mu-1}(z\,y),
\end{equation}
where $-1<\Re(\mu)<\Re(2\nu-1/2)$.
\end{cor}
\begin{pf}
We set 
\begin{equation}
\label{besselkMuNupone}
g(x)=x^{\mu}\,J_{\mu}(z\,x)
\end{equation}
in \eqref{l1:1}. Using the relation \eqref{l2:l} and then the formula \cite[Entry (30), p.185]{E1}, we find
\begin{align}
\notag
\mathcal{L}_{2}\big\{x^{\mu}\,J_{\mu}(z\,x);u\big\}
&=
\frac{1}{2}\,\mathcal{L}\big\{x^{\mu/2}\,J_{\mu}(z\,x^{1/2});u^{2}\big\}
\\
\label{besselkMuNuptwo}
&=
\frac{z^{\mu}}{2^{\mu+1}}\,u^{-2\mu-2}\,\exp\bigg(-\frac{z^{2}}{4u^{2}}\bigg),
\end{align}
where $\Re(\mu)>-1$.
Multiplying both sides of Eq. \eqref{besselkMuNuptwo} with $u^{2\,\nu-2}$ and then applying the $\mathcal{L}_{2}$-transform, we obtain
\begin{equation}
\label{besselkMuNupthree}
\mathcal{L}_{2}\Big\{u^{2\,\nu-2}\,\mathcal{L}_{2}\big\{x^{\mu}\,J_{\mu}(z\,x);u\big\};y\Big\}
=
\frac{z^{\mu}}{2^{\mu+1}}\,\mathcal{L}_{2}\Bigg\{u^{2\nu-2\mu-4}\,
\exp\bigg(-\frac{z^{2}}{4u^{2}}\bigg);y\Bigg\}
\end{equation}
Once again using the relation \eqref{l2:l} and then the formula \cite[Entry (29), p.146]{E1}, we obtain the assertion \eqref{besselkMuNuone}.
\qed
\end{pf}
\begin{rem}\label{rbesselkMuNu}
If we use the definition \eqref{pn2} of the $\mathcal{P}_{\nu,2}$-transform, we may write the formula \eqref{besselkMuNuone} of Corollary \ref{besselkMuNu} as
\begin{equation}
\label{rbesselkMuNuone}
\int_{0}^{\infty}
\frac{x^{\mu+1}\,J_{\mu}(z\,x)}{\big(x^{2}+y^{2}\big)^{\nu}}\,dx
=
\frac{1}{\Gamma(\nu)}\,\Big(\frac{z}{2}\Big)^{\nu-1}\,y^{\mu-\nu+1}\,K_{\nu-\mu-1}(z\,y),
\end{equation}
where $-1<\Re(\mu)<\Re(2\nu-1/2)$, (cf. \cite[Entry 6.565 (4),  p. 686]{GR}). 
\end{rem}
\begin{rem}\label{rbesselJNu}
If we put $\nu=\mu+3/2$ in \eqref{rbesselkMuNuone} and use the formula
\begin{equation}
\label{rbesselJNuone}
K_{1/2}(x)=K_{-1/2}(x)=\Big(\frac{\pi}{2x}\Big)^{1/2}\,\exp(-x)
\end{equation}
we obtain
 \begin{equation}
\label{rbesselJNutwo}
\int_{0}^{\infty}
\frac{x^{\mu+1}\,J_{\mu}(z\,x)}{\big(x^{2}+y^{2}\big)^{\mu+3/2}}\,dx=
\frac{\sqrt{\pi}\,z^{\mu}\,\exp(-z\,y)}{2^{\mu+1}\,y\,\Gamma(\mu+3/2)},
\end{equation}
where $\Re(\mu)>-1$, (cf. \cite[Entry 6.565 (3),  p. 686]{GR}). Similarly, setting $\nu=\mu+1/2$ in \eqref{rbesselkMuNuone} and using the formula \eqref{rbesselJNuone}, we obtain 
\begin{equation}
\label{rbesselJNuthree}
\int_{0}^{\infty}
\frac{x^{\mu+1}\,J_{\mu}(z\,x)}{\big(x^{2}+y^{2}\big)^{\mu+1/2}}\,dx
=
\frac{\sqrt{\pi}\,z^{\mu-1}\,\exp(-z\,y)}{2^{\mu}\,\Gamma(\mu+1/2)},
\end{equation}
where $\Re(\mu)>-1/2$, (cf. \cite[Entry 6.565 (2),  p. 686]{GR}).
\end{rem}

\begin{cor}\label{c1}
We have
\begin{equation}
\label{c1:1}
\mathcal{G}\big\{x^{\mu-1}\,J_{\mu}(z\,x);y\big\}
=
\Big(\frac{2}{\pi\,z}\Big)^{1/2}\,y^{\mu+1/2}\,K_{\mu+1/2}(z\,y),
\end{equation}
where $-1<\Re(\mu)<1/2$ and
\begin{equation}
\label{c1:2}
\mathcal{P}\big\{x^{\mu}\,J_{\mu}(z\,x);y\big\}
=
y^{\mu}\,K_{\mu}(z\,y),
\end{equation}
where $-1<\Re(\mu)<3/2$.
\end{cor}
\begin{pf}
Setting $\nu=1/2$ in \eqref{besselkMuNuone} of our Corollary \ref{besselkMuNu} and using the relationship \eqref{pn2:g}, we obtain the special case \eqref{c1:1}. Similarly, the special case \eqref{c1:2} follows upon setting $\nu=1$ in \eqref{besselkMuNuone}, using the relationship \eqref{pn2:w}, and making use of the fact that the function $K_{\nu}(x)$ is an even function with respect to the index $\nu$.
\end{pf}
\begin{thm}\label{t1} 
If the conditions stated in Lemma \ref{l1} are satisfied, then the Parseval-Goldstein type relations
\begin{align}
\label{t1:1}
\int_{0}^{\infty}
y^{2\nu-1}\,\mathcal{L}_{2}\left\{f(x);y\right\}\,\mathcal{L}_{2}\left\{g(u);y\right\}\,dy
&=
\frac{\Gamma(\nu)}{2}
\int_{0}^{\infty}
x\,f(x)\,
\mathcal{P}_{\nu,2}\big\{g(u);x\big\}\,dx
\\
\label{t1:2}
\int_{0}^{\infty}
y^{2\nu-1}\,\mathcal{L}_{2}\left\{f(x);y\right\}\,\mathcal{L}_{2}\left\{g(u);y\right\}\,dy
&=
\frac{\Gamma(\nu)}{2}
\int_{0}^{\infty}
u\,g(u)\,
\mathcal{P}_{\nu,2}\big\{f(x);u\big\}\,du
\\
\intertext{and}
\label{t1:3}
\int_{0}^{\infty}
x\,f(x)\,
\mathcal{P}_{\nu,2}\big\{g(u);x\big\}\,dx
&=
\int_{0}^{\infty}
u\,g(u)\,
\mathcal{P}_{\nu,2}\big\{f(x);u\big\}\,du
\end{align} 
hold true.
\end{thm}
\begin{pf} 
We only give the proof of \eqref{t1:1}, as the proof of \eqref{t1:2} is similar. Identity \eqref{t1:3} follows from the identities \eqref{t1:1} and \eqref{t1:2}. 

Using the definition \eqref{l2} of the $\mathcal{L}_{2}$-transform, we have
\begin{align}
\notag
\int_{0}^{\infty}
y^{2\nu-1}\,\mathcal{L}_{2}&\big\{f(x);y\big\}\,\mathcal{L}_{2}\big\{g(u);y\big\}\,dy
\\
\label{t1:p1}
&=
\int_{0}^{\infty}
y^{2\nu-1}\,\mathcal{L}_{2}\big\{g(u);y \big\}\,
\left[\int_{0}^{\infty}x\,\exp\big(-x^{2}\,y^{2}\big)\,f(x)\,dx\right]\,dy.
\end{align} 
Changing the order of integration (which is permissible by absolute convergence of the integrals involved) and using the definition \eqref{l2} of the $\mathcal{L}_{2}$-transform once again, we find from \eqref{t1:p1} that
\begin{align}
\notag
\int_{0}^{\infty}
y^{2\nu-1}\,\mathcal{L}_{2}&\big\{f(x);y\big\}\,\mathcal{L}_{2}\big\{g(u);y\big\}\,dy
\\
\notag
&=
\int_{0}^{\infty}
x\,f(x)
\left[\int_{0}^{\infty}y^{2\nu-1}\,\exp\big(-x^{2}\,y^{2}\big)
\mathcal{L}_{2}\big\{g(u);y\big\}\,dy\right]\,dx
\\
\label{t1:p2}
&=
\int_{0}^{\infty}
x\,f(x)\,
\mathcal{L}_{2}\Big\{y^{2\nu-2}\,
\mathcal{L}_{2}\big\{g(u);y\big\};x\Big\}\,dx
\end{align} 
Now the assertion \eqref{t1:1} easily follows from \eqref{t1:p2} and \eqref{l1:1} of the Lemma~\ref{l1}.
\qed
\end{pf}
Setting $\nu=1$ in the identities \eqref{t1:1}, \eqref{t1:2}, and \eqref{t1:3} of our Theorem \ref{t1} and using the relations \eqref{pn2:w}, \eqref{l2l2w:1}, we obtain identities involving the $\mathcal{L}_{2}$-transform and the Widder potential transform contained in (cf. \cite[p. 519, Eq. (2.4)]{YS})
\begin{cor}\label{l2w} 
If the conditions stated in Lemma \ref{l1} are satisfied, then the Parseval-Goldstein type relations
\begin{align}
\label{l2w:1}
\int_{0}^{\infty}
y\,\mathcal{L}_{2}\left\{f(x);y\right\}\,\mathcal{L}_{2}\left\{g(u);y\right\}\,dy
&=
\frac{1}{2}
\int_{0}^{\infty}
x\,f(x)\,
\mathcal{P}\big\{g(u);x\big\}\,dx
\\
\label{l2w:2}
\int_{0}^{\infty}
y\,\mathcal{L}_{2}\left\{f(x);y\right\}\,\mathcal{L}_{2}\left\{g(u);y\right\}\,dy
&=
\frac{1}{2}
\int_{0}^{\infty}
u\,g(u)\,
\mathcal{P}\big\{f(x);u\big\}\,du
\\
\intertext{and}
\label{l2w:3}
\int_{0}^{\infty}
x\,f(x)\,
\mathcal{P}\big\{g(u);x\big\}\,dx
&=
\int_{0}^{\infty}
u\,g(u)\,
\mathcal{P}\big\{f(x);u\big\}\,du
\end{align} 
hold true.
\end{cor}
An immediate consequence of Theorem \ref{t1} is contained in
\begin{cor}\label{ck} 
If the integrals involved converge absolutely, then we have
\begin{align}
\mathcal{L}_{2}\bigg\{y^{2\mu-2\nu}\,\mathcal{L}_{2}\Big\{f(x);\frac{1}{2y}\Big\};z\bigg\}
\label{ck:1}
&=
\frac{z^{\nu-\mu-\frac{3}{2}}}{2^{\mu-\nu+1}}\,\mathcal{K}_{\nu-\mu-1}
\big\{x^{\mu-\nu+\frac{3}{2}}\,f(x);z\big\}
\\
\mathcal{L}_{2}\bigg\{y^{2\mu-2\nu}\,\mathcal{L}_{2}\Big\{f(x);\frac{1}{2y}\Big\};z\bigg\}
\label{ck:2}
&=
\frac{2^{2\nu-\mu-2}}{z^{\mu+\frac{1}{2}}}\,\Gamma(\nu)\,\mathcal{H}_{\mu}
\Big\{u^{\mu+\frac{1}{2}}\,\mathcal{P}_{\nu,2}\big\{f(x);u\big\};z\Big\}
\\
\intertext{and}
\mathcal{H}_{\mu}
\Big\{u^{\mu+\frac{1}{2}}\,\mathcal{P}_{\nu,2}\big\{f(x);u\big\};z\Big\}
\label{ck:3}
&=
\frac{1}{\Gamma(\nu)}\,\Big(\frac{z}{2}\Big)^{\nu-1}\,
\mathcal{K}_{\nu-\mu-1}
\Big\{x^{\mu-\nu+\frac{3}{2}}\,f(x);z\Big\},
\end{align}
where $-1<\Re(\mu)<\Re(2\nu-\frac{1}{2})$ and $\Re(\nu)>0$. 
\end{cor}
\begin{pf} 
We set 
\begin{equation}
\label{ck:p1}
g(u)=u^{\mu}\,J_{\mu}(z\,u)
\end{equation}
in \eqref{t1:1} of Theorem \ref{t1}.
Utilizing \eqref{besselkMuNuptwo}, we have
\begin{equation}
\label{ck:p2}
\mathcal{L}_{2}\big\{g(u);y\big\}
=
\frac{z^{\mu}}{2^{\mu+1}}\,y^{-2\mu-2}\,\exp\bigg(-\frac{z^{2}}{4y^{2}}\bigg).
\end{equation}
Utilizing \eqref{besselkMuNuone}, we have
\begin{equation}
\label{ck:p3}
\mathcal{P}_{\nu,2}\big\{g(u);x\big\}
=
\frac{1}{\Gamma(\nu)}\,\Big(\frac{z}{2}\Big)^{\nu-1}\,x^{\mu-\nu+1}\,K_{\nu-\mu-1}(z\,x).
\end{equation}
Substituting the results  \eqref{ck:p1}, \eqref{ck:p2}, and \eqref{ck:p3} into 
\eqref{t1:1} of Theorem \ref{t1}, we obtain
\begin{align}
\notag
\int_{0}^{\infty}
y^{2\nu-2\mu-3}\,&\exp\bigg(-\frac{z^{2}}{4y^{2}}\bigg)\,\mathcal{L}_{2}\left\{f(x);y\right\}\,dy
\\
\label{ck:p4}
&=
\Big(\frac{z}{2}\Big)^{\nu-\mu-1}
\int_{0}^{\infty}
x^{\mu-\nu+2}\,K_{\nu-\mu-1}(z\,x)\,f(x)\,dx.
\end{align}
The assertion \eqref{ck:1} follows if we change the variable of the integration to $y=1/2v$ on the left-hand side of \eqref{ck:p4}; and then use the definition \eqref{l2} of the $\mathcal{L}_{2}$-transform on the left-hand side of \eqref{ck:p4} and the definition \eqref{kn} of the $\mathcal{K}$-transform on the right hand side of \eqref{ck:p4}.

To prove the identity \eqref{ck:2}, we substitute Eqs. \eqref{ck:p1} and \eqref{ck:p2} into  \eqref{t1:2} of Theorem \ref{t1}, we obtain
\begin{align}
\notag
\int_{0}^{\infty}
y^{2\nu-2\mu-2}\,&\exp\bigg(-\frac{z^{2}}{4y^{2}}\bigg)\,\mathcal{L}_{2}\left\{f(x);y\right\}\,dy
\\
\label{ck:p5}
&=
\Big(\frac{2}{z}\Big)^{\mu}\,\Gamma(\nu)\,
\int_{0}^{\infty}
u^{\mu+1}\,J_{\mu}(z\,u)\,\mathcal{P}_{\nu,2}\big\{f(x);u\big\}\,du.
\end{align}
The assertion \eqref{ck:2} follows if we change the variable of the integration to $y=1/2v$ on the left-hand side of \eqref{ck:p5}, then use the definition \eqref{l2} of the $\mathcal{L}_{2}$-transform on the left-hand side of \eqref{ck:p5} and use the definition \eqref{hankel} of the Hankel transform on the right hand side of \eqref{ck:p5}.

The proof of identity \eqref{ck:3} immediately follows from the identities \eqref{ck:1} and \eqref{ck:2}.
\qed
\end{pf}

\begin{rem}\label{rck}
Setting $\nu=1$ in Corollary \ref{ck}, making use of the fact that $K_{\nu}(x)$ is an even function with respect to the index $\nu$ and the relationship \eqref{pn2:w}, we obtain the following identities:
\begin{align}
\mathcal{L}_{2}\bigg\{y^{2\mu-2}\,\mathcal{L}_{2}\Big\{f(x);\frac{1}{2y}\Big\};z\bigg\}
\label{rck:1}
&=
\frac{z^{-\mu-\frac{1}{2}}}{2^{\mu}}\,\mathcal{K}_{\mu}
\big\{x^{\mu+\frac{1}{2}}\,f(x);z\big\}
\\
\mathcal{L}_{2}\bigg\{y^{2\mu-2}\,\mathcal{L}_{2}\Big\{f(x);\frac{1}{2y}\Big\};z\bigg\}
\label{rck:2}
&=
\frac{2^{-\mu}}{z^{\mu+\frac{1}{2}}}\,\mathcal{H}_{\mu}
\Big\{u^{\mu+\frac{1}{2}}\,\mathcal{P}\big\{f(x);u\big\};z\Big\}
\\
\intertext{and}
\mathcal{H}_{\mu}
\Big\{u^{\mu+\frac{1}{2}}\,\mathcal{P}\big\{f(x);u\big\};z\Big\}
\label{rck:3}
&=
\mathcal{K}_{\mu}
\Big\{x^{\mu+\frac{1}{2}}\,f(x);z\Big\},
\end{align}
where $-1<\Re(\mu)<3/2$. We would like to note that the identity \eqref{rck:1} is obtained earlier (see, \cite[p. 519, Eq. (2.11)]{YS}). If we put $\mu=-1/2$ in \eqref{rck:3} and use the special cases \eqref{rbesselJNuone} and 
\begin{equation}
\label{rck:4}
J_{-\frac{1}{2}}(x)=\Big(\frac{2}{\pi\,x}\Big)^{1/2}\,\cos(x).
\end{equation}
we obtain another known identity \cite[p. 153, Eq. 10]{Y91}.
\end{rem}
\begin{cor}\label{cali} 
If the integrals involved converge absolutely, then we have
\begin{align}
\label{cali:1}
\int_{0}^{\infty}y^{2\nu-2\mu-1}\mathcal{L}_{2}\big\{g(u);y\big\}\,dy
&=
\frac{\Gamma(\nu)}{\Gamma(\mu)}\int_{0}^{\infty}x^{2\mu-1}\mathcal{P}_{\nu,2}\big\{g(u);x\big\}\,dx
\\
\label{cali:2}
\int_{0}^{\infty}y^{2\nu-2\mu-1}\mathcal{L}_{2}\big\{g(u);y\big\}\,dy
&=
\frac{\Gamma(\nu-\mu)}{2}\,
\int_{0}^{\infty}u^{2\mu-2\nu-1}\,g(u)\,du
\\
\label{cali:3}
\int_{0}^{\infty}x^{2\mu-1}\,\mathcal{P}_{\nu,2}\big\{g(u);x\big\}\,dx
&=
\frac{1}{2}\,B(\mu,\nu-\mu)\,
\int_{0}^{\infty}u^{2\mu-2\nu-1}\,g(u)\,du
\end{align}
where $0<\Re(\mu)<\Re(\nu)$, and $B(x,y)$ denotes the beta function. 
\end{cor}
\begin{pf} 
The proof of the identity \eqref{cali:1} follows upon setting 
\begin{equation}
\label{calip:1}
f(x)=x^{2\mu-2}
\end{equation}
in \eqref{t1:1} of Theorem \ref{t1} and using the formula \eqref{l1:p3}.

Next we verify  the identity \eqref{cali:2}. We replace $f(x)$ in the assertion \eqref{t1:2} of the Theorem \ref{t1} with the function considered in \eqref{calip:1}. Using the known formula \cite[p. 310, Entry (19)]{E1}, we evaluate the $\mathcal{P}_{\nu,2}$-transform
\begin{equation}
\label{calip:2}
\mathcal{P}_{\nu,2}\big\{x^{2\mu-2};u\big\}=\int_{0}^{\infty}\frac{x^{2\mu-1}}{(x^{2}+u^{2})^{\nu}}\,dx=\frac{1}{2}\,u^{2\mu-2\nu}\,B(\mu,\nu-\mu).
\end{equation}
We have the well known relationship 
\begin{equation}
\label{calip:3}
B(x,y)=\frac{\Gamma(x)\,\Gamma(y)}{\Gamma(x+y)}
\end{equation}
between the gamma function and the beta function. We substitute the equations \eqref{calip:1} and \eqref{calip:2} into the identity \eqref{t1:2} of the Theorem \ref{t1}. The assertion stated in \eqref{cali:2} of our corollary immediately follows upon using the results \eqref{l1:p3} and \eqref{calip:3} after the substitution.

The last identity \eqref{cali:2} is obtained by using the result \eqref{calip:3} in the identity \eqref{t1:3} of the Theorem \ref{t1}. \qed
\end{pf}

\begin{cor}\label{ci} 
If the integrals involved converge absolutely, then we have
\begin{align}
\label{ci:1}
\notag
\mathcal{P}_{\mu,2}\Big\{\mathcal{P}_{\nu,2}&\big\{g(u);x\big\};t\Big\}\\
&=\frac{1}{\Gamma(\nu)}\int_{0}^{\infty} y^{2\nu+2\mu-2}\,\exp\big(t^{2}\,y^{2}\big)\,\Gamma\big(-\mu+1;t^{2}\,y^{2}\big)\,\mathcal{L}_{2}\big\{g(u);y\big\}\,dy
\\
\label{ci:2}
\notag
\mathcal{G}\Big\{\mathcal{P}_{\nu,2}&\big\{g(u);x\big\};t\Big\}\\
&=\frac{\sqrt{\pi}}{\Gamma(\nu)}\int_{0}^{\infty} y^{2\nu-1}\,\exp\big(t^{2}\,y^{2}\big)\,\text{\rm Erfc}(ty)\,\mathcal{L}_{2}\big\{g(u);y\big\}\,dy
\end{align}
where $0<\Re(\mu)<\Re(\nu)$, and $B(x,y)$ denotes the beta function. 
\end{cor}

\section{Illustrative Examples}
\numberwithin{equation}{section}

An interesting illustration for the identity \eqref{l1:1} asserted by Lemma \ref{l1} 
for the exponential integral $\text{Ei}(x)=-E_{1}(-x)$ defined by
\begin{equation}
\label{e1}
\text{E}_{1}(x)=\int_{x}^{\infty}\frac{e^{-u}}{u}\,du
\end{equation}
is contained in
the following example. 
\begin{exmp} 
We show that
\begin{equation}
\label{ex1:1}
P_{\nu,2}\bigg\{\text{E}_{1}\bigg(\frac{a^{2}}{x^{2}}\bigg);y\bigg\}=
\frac{\Gamma(\nu-1)}{2a(\nu-1)}
\,y^{-2\nu+3}\,\exp\bigg(\frac{a^{2}}{2y^{2}}\bigg)\,W_{-\nu+\frac{3}{2},0}\bigg(\frac{a^{2}}{y^{2}}\bigg),
\end{equation}
where $\Re(\nu)>1$ and $W_{\lambda,\mu}(x)$ denotes Whittaker's function.
\end{exmp}
{\noindent \it Demonstration.} We put
\begin{equation}
\label{ex1:p1} 
g(x)=\text{Ei}\big(-a^{2}/x^{2}\big)
\end{equation}
in the identity \eqref{l1:1} of Lemma \ref{l1}. Using the relation \eqref{l2:l} and the known identity \cite[p. 136, Entry 3.4.1-(13)]{PBM}, we find that
\begin{equation}
\label{ex1:p2} 
\mathcal{L}_{2}\big\{g(x);u\big\}=\frac{1}{2}\,\mathcal{L}\bigg\{\text{Ei}\bigg(-\frac{a^{2}}{x}\bigg);u^{2}\bigg\}=-\frac{1}{u^{2}}\,K_{0}(2\,a\,u).
\end{equation}
Multiplying both sides of \eqref{ex1:p2} by $u^{2\nu-2}$, applying the $\mathcal{L}_{2}$-transform, and then using the relation \eqref{l2:l} once more we deduce that
\begin{equation}
\label{ex1:p2} 
\mathcal{L}_{2}\Big\{\mathcal{L}_{2}\big\{g(x);u\big\};y\Big\}
=-\mathcal{L}_{2}\big\{u^{2\nu-4}\,K_{0}(2\,a\,u);y\big\}
=-\frac{1}{2}\,\mathcal{L}\big\{u^{\nu-2}\,K_{0}(2\,a\,u^{1/2});y^{2}\big\}.
\end{equation}
Now the assertion \eqref{ex1:1} follows from \eqref{ex1:p2} upon using the known formula \cite[p. 353, Entry 3.16.2-(3)]{PBM}.

Another illustration for the identity \eqref{l1:1} asserted by Lemma \ref{l1} 
for the Tricomi confluent hypergeometric function $\Psi(a,b;z)$ is contained in
the following example. The Tricomi hypergeometric function is defined by
\begin{equation}
\label{tricomi}
\Psi(a,b;z)=\frac{\Gamma(1-b)}{\Gamma(a-b+1)}\,{}_{1}F_{1}(a,b;z)+
\frac{\Gamma(b-1)}{\Gamma(a)}\,z^{1-b}\,{}_{1}F_{1}(a-b+1,2-b;z),
\end{equation}
where $b\notin\mathbb{Z}$ and ${}_{1}F_{1}(a,b;z)$ is the Kummer confluent hypergeometric function defined by 
\begin{equation}
\label{kummer}
{}_{1}F_{1}(a,b;z)=\sum_{k=0}^{\infty}\frac{(a)_{k}\,z^{k}}{(b)_{k}\,k!}.
\end{equation}

\begin{exmp} 
We show that
\begin{equation}
\label{extr:1}
P_{\nu,2}\bigg\{x^{2\,\mu}\,\exp\big(-a^{2}\,x^{2}\big);y\bigg\}=
\frac{\Gamma(\mu+1)}{2}
\,a^{2\nu-2\mu-2}\,\Psi\big(\nu,\nu-\mu;a^{2}\,y^{2}\big),
\end{equation}
where $\Re(\mu)>-1$, $\Re(\nu)>0$, and $\nu-\mu\notin\mathbb{Z}$.
\end{exmp}
{\noindent \it Demonstration.} If we set
\begin{equation}
\label{extr:d1}
g(x)=x^{2\mu}\,\exp\big(-a^{2}\,x^{2}\big)
\end{equation}
in the assertion \eqref{l1:1} of Lemma \ref{l1} and use the relationship \eqref{l2:l} and the known formula \cite[p. 144, Entry (3)]{E1}, we find that
\begin{equation}
\label{extr:d2}
\mathcal{L}_{2}\big\{g(x);u\big\}=\frac{1}{2}\,\mathcal{L}\big\{x^{\mu}\,\exp\big(-a^{2}\,x\big);u^{2}\big\}
=\frac{1}{2}\,\Gamma(\mu+1)\,\big(u^{2}+a^{2}\big)^{-\mu-1}.
\end{equation}
Multiplying both sides of \eqref{extr:d2} by $u^{2\nu-2}$, applying the $\mathcal{L}_{2}$-transform and using the relationship \eqref{l2:l} once more, we deduce that
\begin{align}
\notag
\mathcal{L}_{2}\bigg\{u^{2\nu-2}\mathcal{L}_{2}\big\{g(x);u\big\};y\bigg\}
&=\frac{\Gamma(\mu+1)}{2}\,\mathcal{L}_{2}\Big\{u^{2\nu-2}\,\big(u^{2}+a^{2}\big)^{-\mu-1};y\Big\}
\\
\label{extr:d3}
&=\frac{\Gamma(\mu+1)}{4}\,\mathcal{L}\Big\{u^{\nu-1}\,\big(u+a^{2}\big)^{-\mu-1};y^{2}\Big\}.
\end{align}
Now the assertion \eqref{extr:1} follows upon using the known formula \cite[p. 18, Entry 2.1.3-(1)]{PBM} to evaluate the Laplace transform on the right hand side of \eqref{extr:d3}.

\begin{exmp} 
We show that
\begin{align}
\label{exk:1}
\mathcal{K}_{\nu-\mu-1}\Big\{x^{\mu-\nu+\frac{1}{2}}\,\cos(a\,x);z\Big\}&=
\frac{2^{\mu-\nu}\,\sqrt{\pi}}{z^{\nu-\mu-\frac{3}{2}}}
\,
\frac{\Gamma\big(\mu-\nu+\frac{3}{2}\big)}{\big(z^{2}+a^{2}\big)^{\mu-\nu+\frac{3}{2}}}
\\
\label{exk:2}
\mathcal{H}_{\mu}\Big\{u^{\mu+\nu}\,K_{\nu-\frac{1}{2}}(a\,u);z\Big\}
&=
2^{1/2}\Big(\frac{a}{2}\Big)^{\nu-\frac{1}{2}}\,(2z)^{\mu+\frac{1}{2}}
\,
\frac{\Gamma\big(\mu-\nu+\frac{3}{2}\big)}{\big(z^{2}+a^{2}\big)^{\mu-\nu+\frac{3}{2}}},
\end{align}
provided that the conditions of Corollary \ref{ck} hold true.
\end{exmp}
{\noindent \it Demonstration.} If we put 
\begin{equation}
\label{exk:d1}
f(x)=\frac{1}{x}\,\cos(a\,x),
\end{equation}
in our Corollary \ref{ck}, then using the relationship \eqref{l2:l} and the known formula \cite[p. 158, Entry 4.7 (67)]{E1} we find that
\begin{equation}
\label{exk:d2}
\mathcal{L}_{2}\Big\{f(x);\frac{1}{2y}\Big\}
=\frac{1}{2}\,\mathcal{L}\Big\{x^{-1/2}\,\cos\big(a\,x^{1/2}\big);\frac{1}{4y^{2}}\Big\}
=\sqrt{\pi}\,y\,\exp\big(-a^{2}\,y^{2}\big).
\end{equation}
Multiplying both side of \eqref{exk:d2} by $y^{2\mu-2\nu}$, applying both sides the $\mathcal{L}_{2}$-transform, and finally using the relation \eqref{l2:l} we obtain 
\begin{align}
\notag
\mathcal{L}_{2}\bigg\{y^{2\mu-2\nu}\,\mathcal{L}_{2}\Big\{f(x);\frac{1}{2y}\Big\};z\bigg\}
&=\sqrt{\pi}\,\mathcal{L}_{2}\Big\{y^{2\mu-2\nu+1}\,\exp\big(-a^{2}\,y^{2}\big);z\Big\}\\
\notag
&=\frac{\sqrt{\pi}}{2}\,\mathcal{L}\Big\{y^{\mu-\nu+\frac{1}{2}}\,\exp\big(-a^{2}\,y\big);z^{2}\Big\}
\label{exk:d3}
\\
&=\frac{\sqrt{\pi}}{2}\,\frac{\Gamma\Big(\mu-\nu+\frac{3}{2}\Big)}{\Big(z^{2}+a^{2}\big)^{\mu-\nu+\frac{3}{2}}}\cdot
\end{align}
The assertion \eqref{exk:1} follows upon substituting \eqref{exk:d1},  and \eqref{exk:d3} into the identity \eqref{ck:1} of the Corollary \ref{ck}. In order to verify the assertion \eqref{exk:2}, we use the known formula \cite[p. 959, Entry 8.432-5]{GR} to get
\begin{equation}
\label{exk:d4}
\mathcal{P}_{\nu,2}\Big\{\frac{\cos(a\,x)}{x};u\Big\}
=\frac{\sqrt{\pi}}{\Gamma(\nu)}\,\Big(\frac{u}{2a}\Big)^{\nu-\frac{1}{2}}\,K_{\nu-\frac{1}{2}}(a\,u).
\end{equation}
Now the assertion \eqref{exk:2} follows upon substituting \eqref{exk:d1}, \eqref{exk:d3}, and \eqref{exk:d4} into the identity \eqref{ck:2} of the Corollary \ref{ck}.

We conclude this investigation by remarking that many other inÞnite 
integrals can be evaluated in this manner by applying the above Lemma, the 
above Theorem, and their various corollaries and consequences considered 
here.





\begin{thebibliography}{00}


\bibitem{YS} 
O. Y\"urekli and I. Sadek,
A Parseval-Goldstein type theorem on the Widder potential transform and its applications, 
{\it Internat. J. Math and Math. Sci.} 
{\bf 14} (1991), 517--524.





\bibitem{BDY} D. Brown, N. Dernek and O. Y{\"u}rekli, 
Identities for the $E_{2,1}$-transform and their applications,
{\it Appl. Math. Comput.} {\bf 187} (2007), 1557-1566.



\bibitem{DSY} N. Dernek, H. M. Srivastava and O. Y{\"u}rekli,
Parseval-Goldstein type identities involving the $\mathcal{L}_4$-transform and $\mathcal{P}_{4}$-transform and their applications,  {\it Integral Transform. Spec. Funct.} {\bf 18} (2007), 245-253.


\bibitem{GL} 
S. Goldstein, 
Operational representations of Whittaker's  confluent  hypergeometric 
function  and Weber's parabolic cylinder function, {\it Proc. London Math. Soc.
(Ser. 2)} {\bf 34} (1932), 
103--125.





\bibitem{Y99a} O. Y\"{u}rekli,
Theorems on $\mathcal{L}_{2}$-transforms and its applications, 
{\it Complex Variables Theory Appl.}
{\bf 38} (1999), 95--107.



\bibitem{Y99b} 
O. Y\"{u}rekli,
 New identities involving the Laplace and the
$\mathcal{L}_{2}$-transforms and their applications,  
{\it Appl. Math. Comput.} 
{\bf 99}  (1999), 141--151.






\bibitem{WY02} 
O. Y\"urekli and S. Wilson, 
A new method of solving Bessel's differential equation using the $\mathcal{L}_2$-transform.
{\it Appl. Math. Comput.}
{\bf 130} (2002), 587--591.



\bibitem{WY03}
O. Y\"urekli and S. Wilson,
A new method of solving Hermite's differential equation using the $\mathcal{L}_2$-transform.
{\it Appl. Math. Comput.}
{\bf 145} (2003), 495--500.

\bibitem{SY95} 
H. M. Srivastava and O. Y\"urekli,
A theorem on a Stieltjes-type
integral transform and its applications,
 {\it Complex Variables Theory Appl.} 
 {\bf 28} (1995), 159--168.



\bibitem{Y89}
O. Y\"urekli A Parseval-type theorem applied to certain integral transforms, {\it IMA J. Appl. Math.} {\bf 42} (1989),
241--249.



\bibitem{YG} 
O. Yurekli and C. Graziadio,
A theorem on the Laplacae transform and its applications, 
{\it Internat. J. Math. Ed. Sci. Tech.} {\bf 28} (1997), 616--621.

\bibitem{YS98} 
O. Y\"urekli and \"O. Sayg{\i}nsoy,
A theorem  on a Laplace-type integral
transform and its applications. 
{\it Internat. J. Math. Ed. Sci. Tech.} 
{\bf 29} (1998), 561--567.



\bibitem{WD} 
D. V. Widder, 
A transform related to the Poisson integral for a half-plane, {\it Duke Math.
J.} {\bf 33} (1966), 
355--362.


\bibitem{HW} 
I. I. Hirschman and D. V. Widder, \textit{The Convolution Transform}, Princeton University Press, Princeton, New Jersey, 1955.

\bibitem{Gl}
M. L. Glasser, Some Bessel function integrals, {\it Kyungpook Math. J.} {\bf 13} (1973), 171--174.

\bibitem{KSY}
Y. Kahramaner, H. M. Srivastava and O. Y\"urekli,
A theorem on the Glasser transform and its applications.
 {\it Complex Variables Theory Appl.}
{\bf 27} (1995),  7--15.

\bibitem{E1} 
A. Erd\'elyi, W. Magnus, F. Oberhettinger and F. G. Tricomi, 
{\it Tables of Integral Transforms}, Vol. I,
McGraw-Hill Book Company, New York, Toronto and London, 1954.

\bibitem{GR} 
I. S. Gradshteyn, I. M. Ryzhik, 
{\it Table of Integrals, Series, and Products}, 5th edition, 
Academic Press, Inc., Orlando, San Diego, New York, London, Toronto, Montreal, Sydney and Tokyo, 1980.


\bibitem{PBM} 
A.P. Prudnikov, Yu. A. Bryckov and D. I. Marichev,
{\it Integrals and Series}, Academic Press, Vol. 4, 
Gordon and Breach Science Publishers, New York, Reading, Paris, Montreux, Tokyo and Melbourne, 1992.

\bibitem{Y91} 
O. Y\"urekli, An inversion theorem for the Widder potential transform, {\it Internat. J. Math. Ed. Sci. Tech.} {\bf 23} (1992), 152--154.
\end{thebibliography}
\end{document}